\documentclass{elsart}
\usepackage{amssymb}
\usepackage{epsfig}
\usepackage{graphicx}
\journal{Journal of Algebra}

\newcommand\Z{{\mathbb Z}}
\begin{document}
\begin{frontmatter}

\date{May 26, 2010}
\title{Homology operations on homology of quandles}
\author{M. Niebrzydowski}
\address{University of Louisiana at Lafayette, Department of Mathematics,
	 1403 Johnston St., Lafayette, LA 70504-1010}
\ead{mniebrz@gmail.com}

\author{J.H. Przytycki\corauthref{cor}}
\address{The George Washington University, Department of Mathematics,
         Monroe Hall, Room 240,
	 2115 G St. NW, Washington, D.C. 20052}
\corauth[cor]{Corresponding author.}
\ead{przytyck@gwu.edu}

\begin{abstract}
We introduce the concept of the quandle partial derivatives, and use them 
to define extreme chains that yield homological operations.
We apply this to a large class of finite and infinite quandles to show, in particular, 
that they have nontrivial elements in the third and fourth quandle homology.
\end{abstract}

\end{frontmatter}

Quandles are algebraic structures introduced by David Joyce in 
his 1979 Ph.D. thesis \cite{JoyTh} as a powerful 
tool for classifying knots (compare \cite{Joy} and S. Matveev \cite{Mat}). 
Rack homology and homotopy theory were first defined and 
studied in \cite{FRS}, and a modification to quandle homology theory 
was given in \cite{CJKLS} to define knot invariants in a state-sum form 
(so-called cocycle knot invariants). 
In this paper, we consider various homological operations on homology of quandles. 
We introduce the notion of quandle partial derivatives, and extreme chains on which 
appropriate partial derivatives vanish. We also consider the degree one homology 
operations created using elements of the quandle satisfying the so-called $k$-condition.

\section{Definitions and preliminary facts}  \label{1}

\begin{defn}
{\rm A {\it quandle}, $X$, is a set with a binary operation 
$(a, b) \mapsto a * b$
such that
\begin{enumerate}
\item For any $a \in X$, 
$a* a =a$.

\item For any $a,b \in X$, there is a unique $c \in X$ such that 
$a= c*b$. 

\item For any $a,b,c \in X$,
$ (a*b)*c=(a*c)*(b*c)$ (right distributivity). 
\end{enumerate}}
\end{defn}

Note that the second condition can be replaced with the following requirement:
the operation $*_b\colon Q\to Q$, defined by $*_b(x)=x*b$, is a bijection.

\begin{defn}
{\rm A {\it rack} is a set with a binary operation that satisfies 
(2) and (3).}
\end{defn}

We use a standard convention for products in non-associative algebras,
called the left normed convention, that is, whenever parentheses
are omitted in a product of elements $a_1$, $a_2,\ldots,$ $a_n$ of $Q$
then \[a_1*a_2*\ldots *a_n=((\ldots ((a_1*a_2)*a_3)*\ldots)*a_{n-1})*a_n\]
(left association). For example, $a*b*c=(a*b)*c$.

According to \cite{F-R}, the earliest discussion on racks is in the correspondence between 
J.~H.~Conway and G.~Wraith, 
who studied racks in the context of the conjugacy operation in a group. 
They regarded a rack as the wreckage of a group left behind after the group operation 
is discarded and only the notion of conjugacy remains.

There are several basic constructions of quandles from groups and
$\Z[t^{\pm}]$-modules. We list them below stressing the functorial character of
their definition. In particular, the category of quandles has quandles as objects and
quandle homomorphisms as morphisms.
\begin{defn}\
{\rm \begin{enumerate}
\item[(i)] There is a functor $T$ from the category of  abelian groups  to 
the category of quandles, such that for a group $G$, $T(G)$ is a quandle with 
the same underlying set as $G$ and with quandle operation $*$ given by $a*b= 2b-a$. 
For a group homomorphism $f:G \to H$, we define $T(f)=f$ as a function on the set. 
We check that $T(f)$ is a quandle homomorphism: $T(f)(a*b) = f(a*b)=f(2b-a) =2f(b)-f(a)= 
f(a)*f(b) = T(f)(a)*T(f)(b)$.\\
This construction was first considered by M.~Takasaki \cite{Tak}, 
so we call this functor the {\it Takasaki functor}.

\item[(ii)] There is a functor from the category of groups into the category of quandles
in which $a*b= b^{-1}ab$. This functor is called the {\it conjugacy functor}.

\item[(iii)] There is a functor from the category of groups into the category of quandles, in which $a*b= ba^{-1}b$. This functor, generalizing Takasaki functor,
is called the {\it core functor}.

\item[(iv)] There is a functor from the category of $\Z[t^{\pm 1}]$-modules
into the category of quandles in which $a*b= (1-t)b +ta$. This functor, also
generalizing Takasaki functor, is called the {\it Alexander functor}.\\
In particular, if $f\colon A_1\to A_2$ is a $\Z[t^{\pm 1}]$-homomorphism, then it is also an Alexander quandle homomorphism.

\end{enumerate}}
\end{defn}

We recall below the notion of rack and quandle homology. It is useful, following Kamada \cite{Kam}, 
to place it in a slightly more general setting, in which we deal with a rack 
and a rack-set on which the rack acts:
\begin{defn}\label{Definition 1.4}\
{\rm\begin{enumerate}
\item[(1)] For a rack (or a quandle) $X$, the  set $Y$ is a {\it rack-set} (or {\it $X$-set})  
if there is a map $*: Y\times X \to Y$, such that
\begin{enumerate}
\item[(i)] the map $*_x: Y \to Y$, given by $*_x(y) = y*x$, is a bijection, and
\item[(ii)] $(y*a)*b = (y*b)*(a*b)$.
\end{enumerate}

\item[(2)]
For a given rack $X$ and a rack-set $Y$, let 
$C^R_n(X,Y)$ be the free abelian 
group generated by $n$-tuples $(y,x_2,\ldots,x_n)$, with $y\in Y$ and $x_i\in X$, $i=2,\ldots,n$; 
in other words, $C^R_n(X,Y) = {\Z}(Y\times X^{n-1}) = {\Z}Y \otimes {\Z}X^{\otimes n-1}$.\\ 
Define a boundary  homomorphism $\partial: C^R_n(X,Y) \to C^R_{n-1}(X,Y)$ by:
$$\partial(y,x_2,\ldots,x_n) = 
\sum_{i=2}^n (-1)^i((y,\ldots,x_{i-1},x_{i+1},\ldots, 
x_n)$$ $$- (y*x_i,x_2*x_i,\ldots,x_{i-1}*x_i,x_{i+1},\ldots,x_n)).$$ 
$(C^R_*(X,Y),\partial)$ is called the {\it rack chain complex of the pair $(X,Y)$}.\\
The homology of this chain complex is called the {\it homology of the pair 
$(X,Y)$}.
\item[(3)] Assume that $X$ is a quandle. Then we have a subchain 
complex of late degenerate elements, 
$C^{DD}_n(X,Y) \subset C^R_n(X,Y)$, generated by $n$-tuples $(y,x_2,\ldots,x_n)$ 
with $x_{i+1}=x_i$ for some $i$. The subchain complex $(C^{DD}_n(X),\partial)$ 
is called the {\it late degenerated chain complex of the pair $(X,Y)$}.
The homology of this chain complex is called the {\it late degenerated homology 
of $(X,Y)$}, and the homology of the quotient chain complex $C^{LQ}_n(X,Y)= 
C^R_n(X,Y)/C^{DD}_n(X,Y)$ is called the {\it late quandle homology 
of the pair $(X,Y)$}.
\item[(4)] If $X$ is a quandle and $Y$ is an $X$-invariant subquandle of $X$ 
(that is, $y*x \in Y$ for any $y \in Y$ and $x \in X$),
then we have  a subchain
complex $C^D_n(X,Y) \subset C^R_n(X,Y)$, generated by $n$-tuples $(x_1,\ldots,x_n)$
with $x_1 \in Y$ and $x_{i+1}=x_i$ for some $i\in\{1,\ldots,n-1\}$. The subchain complex $(C^D_n(X,Y),\partial)$
is called the {\it degenerated chain complex of a quandle pair $(X,Y)$.}\\
The quotient chain complex $C^Q_n(X,Y)=C^R_n(X,Y)/C^D_n(X,Y)$ is
called the {\it quandle chain complex of a pair $(X,Y)$.}\\
If $X=Y$ then we deal with classical quandle homology theory.
\end{enumerate}}
\end{defn}

Free part of homology of finite racks or quandles  ($free(H_*(X))$) was computed in
\cite{L-N,E-G} (lower bounds for Betti numbers were given in \cite{CJKS-1}).

\begin{thm}\cite{CJKS-1,L-N,E-G}\label{orbit number}
Let $\mathcal{O}$ denote the set of orbits of a rack $X$ with
respect to the action of $X$ on itself by the right multiplication (we say that a rack is connected if it has one orbit).\\
Then:
\begin{enumerate} 
\item[(i)] $rank H_n^R(X) = |\mathcal{O}|^n$ for a finite rack $X$;\\
\item[(ii)] $rank H_n^Q(X) = |\mathcal{O}|(|\mathcal{O}|-1)^{n-1}$
for a finite quandle $X$;\\
\item[(iii)] $rank H_n^D(X) = |\mathcal{O}|^n - |\mathcal{O}|(|\mathcal{O}|-1)^{n-1}$
for a finite quandle $X$;\\
\item[(iv)] $rank H_n^{DD}(X) = |\mathcal{O}|^n - |\mathcal{O}|^2(|\mathcal{O}|-1)^{n-2}$
for a finite quandle $X$,
where $H_n^{DD}(X)$ denotes the delayed degenerated homology.
\end{enumerate}

\end{thm}

If we assume that $X$ satisfies certain mild conditions (as in Theorem \ref{gf}), then 
Theorem \ref{orbit number} follows from Theorem \ref{gf} inductively.

We should stress here, that the assumption of $X$ being finite is essential, as 
M.~Eisermann demonstrated that fundamental quandles of nontrivial knots (so, infinite and connected quandles) 
have $H_2^Q$  equal to $\Z$ \cite{Eis}. \\
For the dihedral quandle $R_k= T(\Z_n)$ (that is, the set $\{0,1,...,n-1\}$ with $a*b= 2b-a$ modulo $n$), 
we have:  
\begin{cor}\label{Theorem 1.2} 
\begin{displaymath}
free(H^Q_n(R_k))=
\left\{
\begin{array}{ll}
\Z & for \ n=1,\ k\ odd \\
0 &  for \ n>1,\ k\ odd\ (R_k \ is\ connected)\\
\Z^2 &  for \ k\ even\ (R_k \ has \ two \ orbits)
\end{array}
\right.
\end{displaymath}
\end{cor} 

Useful information concerning torsion of homology of 
racks and quandles was obtained in \cite{L-N,Moc,N-P-2}. 
In particular, it was shown that:
\begin{thm}\label{Theorem 1.7}\
\begin{enumerate}
\item[(i)] \cite{L-N}\ 
\[ H_2^R(X) \cong H_2^Q(X) \oplus \Z\mathcal{O}_X,\]
\[ H_3^R(X) \cong H_3^Q(X) \oplus H_2^Q(X) \oplus \Z\mathcal{O}_X^2.\]
In particular, $H_3^R(R_k) \cong H_3^Q(R_k) \oplus \Z$ for $k$ odd.
\item[(ii)] \cite{N-P-2}\ 
$H_3^Q(R_p) = \Z_p$,                             
and  $\Z_p \subset H_4^Q(R_p)$ for $p$ an odd prime.
\end{enumerate}
\end{thm}

Quandle homology was generalized to twisted quandle homology in \cite{CES}, 
and further in \cite{A-G,Oht} using the notion of quandle algebra.
Some of our results apply to these generalizations but we will not address it in this paper.
\section{Homological operations obtained from $k$-elements}
Generalizing  \cite{N-P-2}, we consider the degree one homological operation 
related to the group homomorphism
$h_a\colon C^R_{n}(X) \to C^R_{n+1}(X),$ given by
$h_a(w)=(w,a)$, for any $a\in X$, and $w\in X^n$. 
In general, the map $h_a$  is not a chain map, so we
need to symmetrize it with respect to another map $*_a: C^R_{n}(X) \to
C^R_{n}(X)$ given by $*_a(w)=w*a$, for any $w\in X^n$, or more precisely,
$*_a(x_1,\ldots,x_n) = (x_1,\ldots,x_n)*a= (x_1*a,\ldots,x_n*a)$. For a 
symmetrization we need an element $a\in X$ which satisfies a $k$-condition, 
that is, for any $x\in X$ we have $x*a^k=x$ (a quandle in which every element  
satisfies the $k$-condition is called a $k$-quandle; in particular 
a 2-quandle is an involutive quandle or a kei). \\ 
For a $k$-element $a$, we consider a function 
\[h'_a = h_a + *_ah_a +*_a*_ah_a+\ldots + (*_a)^{k-1}h_a.\]
The basic properties of the maps $*_a,$ $h_a$ and $h'_a$ are described in 
the following proposition:
\begin{prop}\label{Proposition 2.1}\
\begin{enumerate}
\item[(i)] For any rack $X$, and $a\in X$, the map $*_a\colon C_n^R(X)\to C_n^R(X)$ is a chain
map chain homotopic to the identity.
\item[(ii)] If $a$ satisfies the $k$-condition, then
$h'_a\colon C_n^R(X) \to C_{n+1}^R(X)$ is a chain map. Notice that if $X$ is finite, then any $a \in X$
satisfies a $k$-condition for some $k$.
\item[(iii)] If $X$ is a quandle, then $h_a(C^D_{n}(X)) \subset C^D_{n+1}(X)$,
 $*_a(C^D_{n}(X)) \subset C^D_{n}(X)$. 
Therefore, the maps $h_a$, $h'_a: C^Q_{n}(X) \to C^Q_{n+1}(X)$, 
and $*_a\colon C_n^Q(X)\to C_n^Q(X)$ are well defined.
Furthermore, $*_a$ and $h'_a$ are chain maps, and $*_a$ is chain homotopic to identity.
\item[(iv)] If $a$ and $b$ are in the same orbit of $X$, then
$h'_a$ and $h'_b$ induce the same map on homology, that is,
$(h'_a)_*=(h'_b)_*: H^W_{n}(X) \to H^W_{n+1}(X),$ where $W=R$, $Q$ or $D$.
\end{enumerate}
\end{prop}

\begin{pf}(i) $*_a$ is a chain map, because
{\small
\[d(*_a(x_1,\ldots,x_n))= d(x_1*a,\ldots,x_n*a)=\]
\[\sum_{i=2}^k (-1)^i((x_1*a,\ldots,x_{i-1}*a,x_{i+1}*a,\ldots,x_n*a)\]
\[-(x_1*a*(x_i*a),\ldots,x_{i-1}*a*(x_i*a),x_{i+1}*a,\ldots,x_n*a))=\]
\[\sum_{i=2}^k (-1)^i((x_1*a,\ldots,x_{i-1}*a,x_{i+1}*a,\ldots,x_n*a)\]
\[-(x_1*x_i*a,\ldots,x_{i-1}*x_i*a,x_{i+1}*a,\ldots,x_n*a))=\]
\[*_a(\sum_{i=2}^k (-1)^i((x_1,\ldots,x_{i-1},x_{i+1},\ldots,x_n)\] 
\[-(x_1*x_i,\ldots,x_{i-1}*x_i,x_{i+1},\ldots,x_n)))=
*_a(d(x_1,\ldots,x_n)).\]}
The homomorphism $(-1)^{n+1}h_a\colon C^R_{n}(X) \to C^R_{n+1}(X)$ is
a chain homotopy between ${\mbox{\rm Id}}$ and $*_a$ chain maps. Namely:
\[d((-1)^{n+1}h_a(x_1,\ldots,x_n))=d((-1)^{n+1}(x_1,\ldots,x_n,a))=\]
\[(-1)^{n+1}((d(x_1,\ldots,x_n),a) +
(-1)^{n+1}((x_1,\ldots,x_n)- (x_1*a,\ldots,x_n*a)))=\]
\[-(-1)^n h_a(d(x_1,\ldots,x_n)) + (Id-*_a)(x_1,\ldots,x_n),\] as needed. In particular, we have $dh_a=h_ad+(-1)^{n+1}(Id-*_a)$.\\
(ii) $dh'_a =d(h_a +*_ah_a+\ldots+(*_a)^{k-1}h_a)= 
dh_a + (dh_a)*a +\ldots+ (dh_a)*a^{k-1} =h_ad +(-1)^{n+1}(Id -*_a) +
(h_ad +(-1)^{n+1}(Id -*_a))*a+\ldots+(-1)^{n+1}(Id -*_a))*a^{k-1} = 
h_ad + (h_a*a)d +\ldots+ (h_a*a^{k-1})d= h'_ad$.\\
(iii) It follows from the definition of rack, degenerate and quandle
chain complex of $X$.\\
(iv) It suffices to consider the case, when there is $x\in X$,
such that $b=a*x$
(the case $b=a\bar* x$ is equivalent to $a=b*x$ so there is no need to consider it separately).\\
Notice, that \[*_x(h'_a(w))= (h'_a(w))*x=
((w+w*a+...+w*a^{k-1}),a)*x=\] 
\[ ((w*x+(w*a)*x+\ldots+(w*a^{k-1})*x),b)=\]
\[((w*x+(w*x)*b+\ldots+(w*x)*b^{k-1}),b)=\]
\[h'_b(w*x)=h'_b(*_x(w)).\]
On the other hand, if $w$ is a cycle then, by (i), $*_x(w)$
is a homologous cycle to $w$. Therefore, $h'_b(*_x(w))$ and $h'_b(w)$ are
homologous by (ii).
Similarly, $*_x(h'_a(w))$ and $h'_a(w)$ are homologous. Therefore, $h'_b(w)$
and $h'_a(w)$ are homologous.
\end{pf}

\begin{defn}\label{Definition 2.2}
{\rm
(a) Let $x_0$ and $x$ be two elements of a rack $X$. We say that $x_0$ is  {\it $j$-connected} to $x$  
if there are elements $x_1,x_2,\ldots,x_j$ such that $x_0*x_1*\ldots*x_j=x$. Notice that if $X$ is 
a quandle, then if $x_0$ is  $j$-connected to $x$ then it is also $(j+1)$-connected to $x$ 
(i.e., $x_0*x_1*\ldots*x_j*x=x$). 
We also say that $x$ is $0$-connected to itself.\\
(b) We say that a rack $X$ is {\it $j$-connected with respect to $x$} if for every $x_0 \in X$, 
$x_0$ is  $j$ (or less)-connected to $x$. We also define the distance $\rho (x_0,x)$ as 
the minimal $j$ such that $x_0$ is $j$-connected to $x$. If such $j$ does not exist, we write 
$\rho (x_0,x)= \infty$. \\
(c) We can visualize definitions (a) and (b) by constructing the oriented graph $G_X$ (Cayley digraph of 
a rack with all elements as generators), whose vertices 
are elements of $X$, and edges (starting from $x$ and ending at $z$) are triples $(x,y,z)$ satisfying,
$x*y=z$. We allow multiple edges and loops. In particular, if $X$ is a quandle, any vertex $x$ has 
a loop $(x,x,x)$. \\
(d) Recall that $\alpha: X \to X$ is an inner automorphism of length $j$ if there are elements
$x_1,x_2,\ldots,x_j$ such that $\alpha(x)= x*x_1*\ldots*x_j$. 
 We say that $X$ is {\it $j$-regular} at $a\in X$ if for any pair of inner automorphisms of length $j$, 
$\alpha_1$ and $\alpha_2$, such that $\alpha_1(x)=\alpha_2(x)=a$ for some $x$, we have $\alpha_1 = \alpha_2$.
  Notice that if $X$ is a quandle, then 
a $j$-regular inner automorphism at $a$ is also a $(j-1)$-regular inner automorphism at $a$. 
Indeed, if $x_0*x_1*\ldots*x_{j-1}=x_0*y_1*\ldots*y_{j-1}=a$, then 
$x_0*x_1*\ldots*x_{j-1}*a=x_0*y_1*\ldots*y_{j-1}*a=a$,  
and by $j$-regularity, for any $x$ we have $x*x_1*\ldots*x_{j-1}*a=x*y_1*\ldots*y_{j-1}*a=a$.
Finally,  
$x*x_1*\ldots*x_{j-1}=x*y_1*\ldots*y_{j-1}=a$, showing $(j-1)$-regularity.\\ 
Notice also, that when our $j$-regular condition is translated to 
Cayley digraphs, then it can be interpreted as ``if two paths of length $j$ starting at $x_0$ have the
common endpoint, then the corresponding paths at any vertex have the same endpoint."\\
(e) The rack $X$ has the {\it quasigroup property}, if for any $a$, $b\in X$, the equation $a*x=b$ has exactly one solution.
Notice that for a connected finite rack, the quasigroup property is equivalent to $X$ being $1$-regular.}
\end{defn}

\begin{exmp}\
{\rm (1) Alexander quandles are $j$-regular for any $j$. To see this let 
$\alpha_1(x) = x*x_1*\ldots*x_j $ and $\alpha_2(x) = x*y_1*\ldots*y_j$. Then we have 
\[ \alpha_1(x) = (1-t)(x_j+tx_{j-1}+\ldots+ t^{j-1}x_1) + t^jx,\]
\[ \alpha_2(x) = (1-t)(y_j+ty_{j-1}+\ldots+ t^{j-1}y_1) + t^jx.\]  
Thus, $\alpha_1(x)- \alpha_2(x) = (1-t)((x_j - y_j)+t(x_{j-1}- y_{j-1})+\ldots+ t^{j-1}(x_1-y_1)),$ which 
does not depend on $x$, proving $j$-regularity of Alexander quandles.\\
(2) Orbits of an Alexander quandle, $A$, are in bijection with the elements of $A/(1-t)A$. Furthermore, 
 every orbit of an Alexander quandle is 1-connected (and a quasigroup if the orbit is finite).\\
To see this, first notice that  $x*y - x= (1-t)(y-x)$. Thus, elements in the same orbit are equal in $A/(1-t)A$. 
Conversely, if $y-x = (1-t)z$, then $y = x*(z+x)$. Thus $y$ and $x$ are in the same orbit and 
$y$ is $1$-connected to $x$ ($\rho(x,y)=1$), as required. Uniqueness of $u$ with $x*u=y$ follows from the
finiteness of the orbit (for an endomorphism of a finite set, an epimorphism is equivalent to a monomorphism). 
In the infinite case, an orbit is not necessary a quasigroup. This happens if $1-t$ annihilates some 
nonzero element of $A$.} 
\end{exmp}

In the following theorem, we generalize Theorem 9 of \cite{N-P-2}.
\begin{thm}\label{Theorem 2.4}
Let $Q$ be a quandle consisting of  $s$ orbits with fixed points $\mathbf{a}=(a_1,a_2,\ldots,a_s)$, 
with one $a_i$ chosen 
from each orbit and satisfying the $k$-condition ($x*a_i^k=x$ for any $x$). 
We define a map $h'_{\mathbf{a}}\colon C^R_n(Q) \to C^R_{n+1}(Q)$ by 
$h'_{\mathbf{a}} = \sum_{i=1}^s h'_{a_i}$. It is a chain map as a sum of chain maps.
\begin{enumerate} 
\item[(1)] If $Q$ is $j$-regular at each $a_i$, and  each orbit of $Q$ is a quasigroup, 
then there is a chain map $\bar h'_{\mathbf{a}}\colon C^R_{n+1}(Q) \to C^R_{n}(Q)$ such that on homology
$(\bar h'_{\mathbf{a}} h'_{\mathbf{a}})_* = sk^2 Id$.
\item[(2)] If each orbit of $Q$ is $j$-connected to $a_i$, and $Q$ is $(j+1)$-regular at $a_i$, 
then there is a chain map $\bar h'_{\mathbf{a}}\colon C^R_{n+1}(Q) \to C^R_{n}(Q)$ such that 
$(\bar h'_{\mathbf{a}} h'_{\mathbf{a}})_* = sk^2 Id$ on homology.

\item[(3)] In particular, if $H^R_{n}(Q)$ has no element
of order dividing $sk$, then $(h'_{\mathbf{a}})_*: H^R_n(Q) \to H^R_{n+1}(Q)$ is a monomorphism.

\end{enumerate}
\end{thm}

\begin{pf} We prove (2) of which (1) is a special case (that includes Alexander quandles).\\
Define $\bar h_{\mathbf{a}}(x_1,\ldots,x_n,x_{n+1}) = (x_1,\ldots,x_n)*u_1*\ldots*u_j$, where $x_{n+1}*u_1*\ldots*u_j=a_i$ for some $i$. 
The map is well defined because $Q$ is $j$-regular at $a_i$.\\
In this situation, we define \[\bar h'_{\mathbf{a}}(x_1,\ldots,x_n,x_{n+1}) = \sum_{t=0}^{k-1}\bar h_{\mathbf{a}}(x_1,\ldots,x_n,x_{n+1})*a_i^t.\] 
We  prove that   $\bar h'_{\mathbf{a}}: C^R_{n+1}(Q) \to C^R_{n}(Q)$ is a chain map. We have:\\  
\[d\bar h'_{\mathbf{a}}(x_1,\ldots,x_n,x_{n+1}) =d(\sum_{t=0}^{k-1}(x_1,\ldots,x_n)*u_1*\ldots*u_j*a_i^t)=\]\[
\sum_{t=0}^{k-1} d(x_1,\ldots,x_n)*u_1*\ldots*u_j*a_i^t.\]
On the other hand, we have:
\[\bar h'_{\mathbf{a}}d(x_1,\ldots,x_n,x_{n+1})=\] 
\[\bar h'_{\mathbf{a}}(d(x_1,\ldots,x_n),x_{n+1}) + 
(-1)^{n+1}\bar h'_\mathbf{a}((x_1,\ldots,x_n)-(x_1,\ldots,x_n)*x_{n+1})=\]
\[\sum_{t=0}^{k-1}d(x_1,\ldots,x_n)*u_1*\ldots*u_j*a_i^t\]\[ + (-1)^{n+1}\bar h'_{\mathbf{a}}((x_1,\ldots,x_n)-(x_1,\ldots,x_n)*x_{n+1}).\]
Finally, we show that $(j+1)$-regularity can be use to show that the last term is equal to zero. Namely, 
for some $a_{i'}$ depending on $x_n$, there are $y_1,\ldots,y_j$ and $z_1,\ldots,z_j$ such that we have 
$x_n*y_1*\ldots*y_j=a_{i'}$ (and so $x_n*y_1*\ldots*y_j*a_{i'}=a_{i'}$), and $x_n*x_{n+1}*z_1*\ldots*z_j=a_{i'}$.
Thus by $(j+1)$-regularity at $a_{i'}$, we have, for any $x$: 
\[x*y_1*\ldots*y_j*a_{i'} = x*x_{n+1}*z_1*\ldots*z_j.\] Thus, 
\[\bar h'_{\mathbf{a}}((x_1,\ldots,x_n)-(x_1,\ldots,x_n)*x_{n+1})=0,\] 
and $\bar h'_\mathbf{a}$ is a chain map.

Then, the formula  $(\bar h'_\mathbf{a} h'_\mathbf{a})_*= sk^2Id$ on homology is easy to check:
\[\bar h'_{\mathbf{a}} h'_{\mathbf{a}}(x_1,\ldots,x_n)=\] 
\[\bar h'_{\mathbf{a}}(\sum_{j=0}^{k-1}\sum_{i=1}^{s}
(x_1,\ldots,x_n,a_i)*a_i^j) = \sum_{u=0}^{k-1}(\sum_{j=0}^{k-1}\sum_{i=1}^{s}(x_1,\ldots,x_n)*a_i^j)*a_i^u =\]
\[k\sum_{j=0}^{k-1}\sum_{i=1}^{s}(x_1,\ldots,x_n)*a_i^j.\]
This holds for any chain, and if $w$ is a cycle in $C^R_n(Q)$, then it is homologous to $w*a_i^j$.
Therefore, $\bar h'_\mathbf{a} h'_\mathbf{a}(w)= k\sum_{j=0}^{k-1}\sum_{i=1}^{s}w*a^i_j$ is homologous to $sk^2w$.

\end{pf}

\begin{cor}\label{Corollary 2.5}\ 
\begin{enumerate} 
\item [(i)] If  $(h'_{\mathbf{a}})_*$ is a non-trivial map (resp. monomorphism), and $X$ is a homogeneous quandle 
(for any two elements $x$ and $y$ of $X$, there is an automorphism $X \to X$ 
sending $x$ to $y$), then $ (h'_{a})_*$ is a nontrivial map (resp. monomorphism) for any $a\in X$.
\item [(ii)] 
If $A$ is a finite  Alexander quandle, then
$(\bar h'_{\mathbf{a}} h'_{\mathbf{a}})_*= sk^2Id$, where $s= |A/((1-t)A)|$ is the number of orbits of $A$, 
and $k$ is the smallest number such that $(1-t^k)$ annihilates $A$.
\end{enumerate}
\end{cor}
\begin{pf} (i)\ 
If $ (h'_{\mathbf{a}})_*(u)\neq 0$, then $ (h'_{a_i})_*(u)\neq 0$ for some $i$,
and thus, by homogeneity of $X$, $(h'_{a})_*(x)\neq 0$ for some $x\in X$. The part (i) of Corollary \ref{Corollary 2.5} follows.\\ 
Part (ii) follows from (i), because any Alexander quandle is homogeneous. Recall that for any $a,b\in A$, 
the module automorphism $f_{a,b}\colon A \to A$ given by $f_{a,b}(x)= x + b -a$ is a quandle automorphism 
sending $a$ to $b$. Also recall that because $x*a^k= x + (1-t^k)(b-x)$, $A$ is a $k$-quandle if and only if 
$(1-t^k)$ annihilates $A$. 
\end{pf}

\begin{exmp}\label{Example 2.6}\
{\rm\begin{enumerate} 
\item [(i)] 
$S_4 = \Z_2[t]/(1+t+t^2)$ is a connected Alexander  
3-quandle with $4$ elements, and $H^R_n(S_4)$ can have only $2$-torsion. 
By Corollary \ref{Corollary 2.5}, 
$(h'_a)_*\colon H^R_n(S_4) \to H^R_{n+1}(S_4)$ is a monomorphism.
In particular, 
$tor H^R_n(S_4)$ is nontrivial (contains $\Z_2$) for $n\geq 2$; we have $H^R_2(S_4)= \Z \oplus \Z_2$.

\item [(ii)]  We computed that $H^Q_2(S_4)=\Z_2$, $H^Q_3(S_4)=\Z_2\oplus \Z_4$, $H^Q_4(S_4)=\Z_2^{2}\oplus \Z_4$, 
$H^Q_5(S_4)=\Z_2^{5}\oplus \Z_4$,  $H^Q_6(S_4)=\Z_2^{9}\oplus \Z_4^{2}$, and 
$H^Q_7(S_4)=\Z_2^{17}\oplus \Z_4^{3}$. Generally, we conjecture that 
$tor H^Q_n(S_4)=\Z_2^{f'_n}\oplus \Z_4^{f_n}$, where 
$\{f_n\}$ are ``delayed" Fibonacci numbers,
that is, $f_n = f_{n-1} + f_{n-3}, \ \textrm{and}\ f(1)=f(2)=0,\ f(3)=1;$ compare \cite{N-P-2}.
Furthermore, $g_n=f'_n +2f_n= log_2(|tor H_n^Q(S_4)|)$
satisfies $g_n=g_{n-1}+g_{n-2}+g_{n-4}$, with $g_1=0, g_2=1,g_3=3$, and $g_4=4.$ 
According to this conjecture, we would have 
$H_8^Q(S_4)= \Z_2^{32} \oplus \Z_4^4$.

\item [(iii)] Consider an Alexander quandle $A_{m,p(t)}= \Z_m[t^{\pm 1}]/(p(t))$, 
where $p(t)$ is a polynomial in variable $t$
with the coefficients of the highest and lowest degree terms of $t$ are invertible in $\Z_m$. 
$A_{m,p(t)}$ has $m^{deg (p(t))}$ elements, and $s= gcd(m,p(1))$ orbits. Then, for a connected quandle $A_{m,p(t)}$, 
the map $(h'_a)_*\colon H^R_n(A_{m,p(t)}) \to H^R_{n+1}(A_{m,p(t)})$ is a monomorphism.
In particular, this holds for $p(x)=[k]_t= 1+t+\ldots +t^{k-1}$, with $gcd(m,k)=1$. Here,
$A_{m,[k]_t}$ is a connected $k$-quandle (notice that $t^k-1=(t-1)[k]_t$ annihilates $A_{m,[k]_t}$).
\end{enumerate}}
\end{exmp}

\section{Annihilation of rack homology}

In this section, we offer an improvement of the results presented in  \cite{L-N,N-P-2},
and we use it to show, in particular, that 
$k^{n-2}$ annihilates $tor H_n(R_{2k})$ for $k$ odd (Corollary \ref{Corollary 3.6}).
We conjecture, however, a much stronger result, at least for $k$ prime:

\begin{conj}\label{Conjecture 3.1}\ 
The number $k$ annihilates $tor H_n(R_{2k})$, unless $k=2^t$, $t>1$. \\
The number $2k$ is the smallest number annihilating $tor H_n(R_{2k})$ for $k=2^t$, $t>1$.
\end{conj}
We checked that  $H_3^Q(R_8)= Z^2 \oplus \Z_8^2$, $H_4^Q(R_8)= \Z^2 \oplus \Z_2^4 \oplus \Z_4^4 \oplus \Z_8^2$, and
$H^Q_3(R_{16})= \Z^2 \oplus \Z_{16}^2 $.

To formulate our further results, 
we use definition of an $X$-set $Y$ and the homology of the pair $(X,Y)$ introduced in 
Definition \ref{Definition 1.4}. 

We start from a fact, allowing as to detect 
torsion in $H_4^Q$ of some quandles $R_{4n}$.

Consider two maps $f,g$ defined as follows (\cite{CJKS-2,N-P-2}:\\
$f\colon C_n^R(X,Y) \to C_{n-1}^R(X)$ given by $f(y,x_2,\ldots,x_n)=(x_2,\ldots,x_n)$,
and\\
$g\colon C_{n-1}^R(X) \to C_n^R(X,Y)$  given by
$g(x_2,\ldots,x_n)= \sum_{y\in Y}(y,x_2,\ldots,x_n)$.
\begin{prop}\label{Proposition 3.2}\
\begin{enumerate}
\item[(i)] The maps $(-1)^nf$, $(-1)^ng$, $fg$, and $gf$ are chain maps. 
Furthermore, $g$ sends degenerate elements in $C_{n-1}^R(X)$ to late degenerate 
elements in $C_n^R(X,Y)$.  \\
\item[(ii)]  $fg = |Y|Id$ on $C_{n-1}^R(X)$  and 
$gf(y,x_2,\ldots,x_n)=\sum_{y\in Y}(y,x_2,\ldots,x_n)$.
\end{enumerate}
\end{prop}
\begin{pf} (i) We easily check that $df+fd=0$ and $dg+gd = 0$. Thus, $(-1)^nf$, $(-1)^ng$, $fg$,
and $gf$ are chain maps.\\
(ii) It follows from a direct computation. 
\end{pf}
\begin{cor}\label{Corollary 3.3}
If $H^Q_n(R_{4k})$ has a $\Z_{4k}$-torsion, then $H^Q_{n+1}(R_{4k})$ contains $\Z_2$ 
(compare Conjecture \ref{Conjecture 3.1}).
\end{cor}
\begin{pf}
Let $X= R_{4k}$ and $Y=R^{even}_{4k}$ be an orbit composed of even numbers.
Then $|Y|=2k$ so if $u$ is a generator of $\Z_{4k}$ in $H^R_n(R_{4k})$, then, by 
Corollary \ref{Corollary 3.3}, $(gf)_*(u)$ is a nonzero element in $H^R_n(R_{4n})$. Therefore, $f_*(u)$ is 
a nontrivial element in $H^R_{n+1}(R_{4k})$. Because $g_*$ is well defined on quandle homology,
$g_*(u)$ is a nontrivial element in 
$H^Q_{n+1}(R_{4k})$ if $u$ is a generator of $\Z_{4k}$ in $H^Q_n(R_{4k})$. \\
\end{pf}

In the next theorem we analyze the chain map $gf\colon C^R_{n}(X,Y) \to C^R_{n}(X,Y)$, 
and show, in particular, that if $Y$ is an orbit of $X$ which is also a quasigroup, then 
$gf$ is chain homotopic to $|Y|Id$ (compare \cite{L-N,CJKS-1}).
\begin{thm} \label{gf}
Let $(X,Y)$ be an $X$-set and let $X_1$ be a finite invariant
subrack of $X$ (e.g., an orbit of $X$). Consider the map $\phi_{X_1}\colon
C_n^R(X,Y) \to C_n^R(X,Y)$ given by $\phi_{X_1}(y,x_2,\ldots,x_n) =
\sum_{x\in X_1}(y*x,x_2,\ldots,x_n)$. Then
\begin{enumerate}
\item[(1)] $\phi_{X_1}$ is a chain map chain homotopic to $|X_1|Id$.
\item[(2)] Let $_y*\colon X_1 \to Y$ be defined by ${_y*}(x)=y*x$, and assume 
that the cardinality of the set of elements of $X_1$ which send $y$ to $y'$, for any pair $y,y'\in Y$, 
is finite and does not depend on $y$ and $y'$. Let us denote this number by $m$ so that $|X_1|=m|Y|$.
Then $mgf=\phi_{X_1}$, which by (1) is chain homotopic to $|X_1|Id$.\\
In particular:
\item[(3)] $f_*$ and $g_*$ are isomorphisms if homology is taken over
any ring in which $|Q_1|$ is invertible, e.g., rational numbers.
\item[(4)] If $_y*$ is always a bijection (e.g., $Y=X_1$, and $X_1$ is a quasigroup), 
then $gf$ is chain homotopic to $|X_1|Id$.
\end{enumerate}
\end{thm}

\begin{pf} Consider a map $H=H_{X,Y,X_1}\colon C_n^R(X,Y) \to C_{n+1}^R(X,Y)$
given by $$H(y,x_2,\ldots,x_n) = (y,\sum_{x\in X_1}x,x_2,\ldots,x_n).$$ It follows that\\
$$dH(y,x_2,\ldots,x_n)=$$
$$|X_1|(y,x_2,\ldots,x_n) - (y*(\sum_{x\in X_1}x),x_2,\ldots,x_n) -
Hd(y,x_2,\ldots,x_n).$$ Therefore, $dH + Hd = |X_1|Id - \phi_{X_1}$. We conclude
that $\phi_{X_1}$ is a chain map chain homotopic to $|X_1|Id$.
Parts (2), (3) and (4) follow straight from our conditions.
\end{pf}

\begin{cor}\label{Corollary 3.5}\
\begin{enumerate}
\item [(i)]
Suppose that $X_0$ is an invariant subrack of $X$, with
$|X_0|$ elements, that is a quasigroup. Then, if
$tor H_n^{R}(X)$ is annihilated by $N$, then
$tor H^R_{n+1}(X,X_0)$ is annihilated by $N|X_0|$.
\item [(ii)] With the notation and assumptions as in Theorem \ref{gf}(2),
we have that if $tor H_n^{R}(X)$ is annihilated by $N$, then
$tor H_{n+1}^{R}(X,Y)$ is annihilated by $N|X_1|/gcd(N,m)$.
\end{enumerate}
\end{cor}
\begin{pf}
(i) Let $a \in tor H_{n+1}^{R}(X,X_0)$. Then, by assumption, $Nf_*(a)=0$
in $H_n^{R}(X)$. Therefore, 
\[N|X_0|a= Ng_*f_*(a)=g_*(Nf_*(a))=g_*(0)=0.\]
(ii) Let $a \in tor H_{n+1}^{R}(X,Y)$. Then $Nf_*(a)=0$
in $H_n^{R}(X)$. Thus, 
\[0=(m/(gcd(N,m)))g_*(Nf_*(a))= (N/(gcd(N,m)(mg_*f_*(a))=\]
\[(N/(gcd(N,m))|X_1|Id,\]
and $tor H_{n+1}^{R}(X,Y)$ is annihilated by $(N/(gcd(N,m)))|X_1|$.
\end{pf}

\begin{cor}\label{Corollary 3.6} For an even dihedral quandle $R_{2k}$, we have:
\begin{enumerate}
\item[(i)]
$tor H^R_n(R_{2k})$ is annihilated by $k^{n-2}$,
for $k$ odd.
\item[(ii)] $tor H_n(R_{2k})$ is annihilated by $2k^{n-1}$,
for an even $k$.
\end{enumerate}
\end{cor}
\begin{pf}
(i)
If $k$ is odd, then $Q_0=Q_{even}$ composed of even numbers is an invariant
quasigroup subquandle in $Q$. Thus, Theorem \ref{gf}(4) applies.
Furthermore, 
\[H^R_n(R_{2k})= H^R_n(R_{2k},Q_0)\oplus H^R_n(R_{2k},Q_{odd}),\]
where $Q_{odd}$ is an invariant subquandle of $R_{2k}$ composed of odd numbers.
Also, \[H^R_n(R_{2k},Q_{even})=H^R_n(R_{2k},Q_{odd}),\] as the quandle
isomorphism $s_+\colon R_{2k} \to R_{2k}$, given by $s_+(i)=i+1$ (modulo $2k$),
sends $Q_{even}$ to $Q_{odd}$, and $Q_{odd}$ to $Q_{even}$.
Finally, by the result of Greene \cite{Gr}, $tor H^R_2(R_{2k}) = 0$ for $k$ odd. The result follows 
by induction on $n$, starting from $n=2$.\\
(ii) We use Theorem \ref{gf}(2), for $Y=Q_{even}$ and $X_1=R_{2k}$. Then,
$m= |X_1|/|Y|= 2$,
and $gcd(m,N)=2$, if $N$ is even, and $gcd(m,N)=1$ if $N$ is odd.
The rest of the proof follows by induction on $n$, starting from $tor H_1(R_{2k})=0$.
\end{pf}

\begin{conj}\label{small}\
\begin{enumerate}
\item[(i)] $tor H_2^Q(R_{4k})=\Z_2^2$.
\item[(ii)] $H_n^Q(R_{2k})$ is annihilated by $2k^{n-2}$, for an even $k$.
\end{enumerate}
\end{conj}
We confirmed Conjecture \ref{small}(i) using GAP \cite{GAP} for $k=1,2,3,4,5,6$.
Part (ii) of Conjecture \ref{small} follows from part (i) by Corollary \ref{Corollary 3.5}(ii);
compare Conjecture \ref{Conjecture 3.1}.

\section{Partial derivatives and homological operations}
By a homological operation of degree $k$ we understand 
any homomorphism $H^R_{*} \to H^R_{*+k}$. A pre-homology operation of degree $k$ is 
a chain map\\ 
$h\colon C^R_{*} \to C^R_{*+k}$.

Below we give the necessary conditions 
for the map $h_w\colon C^R_{*} \to C^R_{*+\ell(w)}$ defined by 
$h_w(y,x_2,\ldots, x_n)=(y,x_2,\ldots,x_n,w)$, where $w\in C^R_{\ell(w)}(Q)$, to be 
a pre-homology operation of degree $\ell(w)$. 
To do this, we first define the partial derivatives 
$\frac{\partial^1 }{\partial q}\colon C^R_n(Q) \to C^R_{n-1}(Q)$.

\begin{defn} 
{\rm We will use the following standard notation. 
\[\partial_i^0(x_1,\ldots,x_n)= (x_1,\ldots,x_{i-1},
x_{i+1},\ldots,x_n),\] 
\[\partial_i^1(x_1,\ldots,x_n)= (x_1*x_i,\ldots,x_{i-1}*x_i,
x_{i+1},\ldots,x_n),\]
and 
\[\partial^0=\sum_{i=1}^n(-1)^i\partial_i^0,\  
\partial^1=\sum_{i=1}^n(-1)^i\partial_i^1,\ \partial =
\partial^0 - \partial^1.\]
\begin{enumerate}
\item[(i)] For any $q\in Q$, we define:\\ 
$\frac{\partial^1}{\partial q}\colon C_n^R(Q) \to C_{n-1}^R(Q)$ by 
$\frac{\partial^1(x_1,\ldots,x_n)}{\partial q}= 
\sum_{i=1}^n(-1)^i\partial_i^1\delta_{x_i,q}$, where the Kronecker delta
$\delta_{x_i,q}=1$ if $x_i=q$, and $0$ otherwise.\\
We have $\partial^1=\sum_{q\in Q}\frac{\partial^1}{\partial q}$.
\item[(ii)] For any $q\in Q$, we define:\\
$\frac{\partial^0}{\partial q}\colon C_n^R(Q) \to C_{n-1}^R(Q)$ by
$\frac{\partial^0(x_1,\ldots,x_n)}{\partial q}=
\sum_{i=1}^n(-1)^i\partial_i^0\delta_{x_i,q}$.\\
We have $\partial^0=\sum_{q\in Q}\frac{\partial^0}{\partial q}$.
\item[(iii)] $\frac{\partial^0}{\partial q}$ sends degenerate 
elements to degenerate elements for any rack $Q$. If $Q$ is a 
quandle, then also $\frac{\partial^1}{\partial q}$ sends degenerate
elements to degenerate elements.
\end{enumerate}}
\end{defn}

\begin{lem} For a quandle $Q$, we have:
\begin{enumerate}
\item[(i)]
$\frac{\partial^0 \partial^0 }{\partial q \partial q}=0$,
\item[(ii)] $\frac{\partial^0\partial^0}{\partial q\partial q'}+
\frac{\partial^0\partial^0}{\partial q'\partial q}=0$.

\item[(iii)] 
$\frac{\partial^1 \partial^1 }{\partial q \partial q}=0$. 

\item[(iv)]
$\partial^0 \frac{\partial^1}{\partial q}+ 
\frac{\partial^1}{\partial q}\partial^0 = 0$.

\item[(v)] $\frac{\partial^0\partial^1}{\partial q\partial q}+
\frac{\partial^1\partial^0}{\partial q\partial q}=0$

\end{enumerate}
\end{lem}

\begin{pf} Parts (i) and (ii) follow from a standard calculation, and no structure on $Q$ is needed, 
it may be any set.\\
(iii) $\frac{\partial^1 \partial^1 }{\partial q \partial q}(x_1,\ldots,x_n)= 
\frac{\partial^1 }{\partial q }( \sum_{i=1}^n (-1)^i(x_1*x_i,\ldots,x_{i-1}*x_i,x_{i+1},\ldots,x_n)\delta_{x_i,q})=$
$\sum_{j<i} (-1)^{i+j}((x_1*x_i)(x_j*x_i),\ldots, (x_{j-1}*x_i)(x_j*x_i),x_{j+1}*x_i,\ldots,$
$x_{i-1}*x_i,x_{i+1},\ldots , x_n)\delta_{x_i,q}\delta_{x_j*q,q} + 
\sum_{j>i} (-1)^{i+j-1}(x_1*x_i*x_j,\ldots,x_{i-1}*x_i*x_j,x_{i+1}*x_j,\ldots ,
x_{j-1}*x_j,x_{j+1},\ldots, x_n)\delta_{x_i,q}\delta_{x_j,q} =0$. 
We use in this calculation all axioms of quandle, in particular $q*q=q$ and ($x*q=q \Rightarrow x=q$). 
In the language of face maps (see the next paragraph for a definition), we can write our property 
$\frac{\partial^1 \partial^1 }{\partial q \partial q}=0$
locally as: 
\[for\ j<i:\ (\partial^1_j\partial^1_i)\delta_{x_i,q}\delta_{x_j*q,q}=
(\partial^1_{i-1}\partial^1_j)\delta_{x_i,q}\delta_{x_j,q}.\]
Parts (iv) and (v):\\
Simplifying the notation from the proof of (iii), let us define the face map
$\frac{\partial^1}{\partial q}{(i)}$ as 
$\partial_i^1 \delta_{x_i,q}$. That is, 
\[\frac{\partial^1}{\partial q}{(i)}(x_1,\ldots,x_n)= (x_1*q,\ldots,x_{i-1}*q,
x_{i+1},\ldots,x_n)\] if $x_i=q$, and $0$ otherwise. Then we need to consider some cases.
For $i<j$:
\[\partial^0_{j-1}\frac{\partial^1}{\partial q}{(i)}(x_1,\ldots,x_n)=\]\[ 
(-1)^{i+j-1}(x_1*q,\ldots,x_{i-1}*q,x_{i+1},\ldots,x_{j-1},x_{j+1},\ldots,x_n)\]
if $x_i=q$, and $0$ otherwise.\\
For $i>j$: 
\[\partial^0_{j}\frac{\partial^1}{\partial q}{(i)}(x_1,\ldots,x_n)=\]\[ 
(-1)^{i+j}(x_1*q,\ldots,x_{j-1}*q,x_{j+1}*q,x_{i-1}*q,x_{i+1},\ldots,x_n)\] 
if $x_i=q$, and $0$ otherwise.\\
Similarly, for $i<j$ we have:
\[\frac{\partial^1}{\partial q}{(i)}\partial^0_{j}(x_1,\ldots,x_n)=\]\[ 
(-1)^{i+j}(x_1*q,\ldots,x_{i-1}*q,x_{i+1},\ldots,x_{j-1},x_{j+1},\ldots,x_n)\]
if $x_i=q$, and $0$ otherwise.\\
For $i>j$:
\[\frac{\partial^1}{\partial q}{(i-1)}\partial^0_{j}(x_1,\ldots,x_n)=\]\[
(-1)^{i+j-1}(x_1*q,\ldots,x_{j-1}*q,x_{j+1}*q,x_{i-1}*q,x_{i+1},\ldots,x_n)\]
if $x_i=q$, and $0$ otherwise.\\
Therefore, for $i<j$:
\[\partial^0_{j-1}\frac{\partial^1}{\partial q}{(i)}+ 
               \frac{\partial^1}{\partial q}_{i}\partial^0_{j} = 0,\] and 
for $i>j$:
\[\partial^0_{j}\frac{\partial^1}{\partial q}{(i)} + 
             \frac{\partial^1}{\partial q}{(i-1)}\partial^0_{j} = 0.\]
Thus, for any rack it follows that 
\[\partial^0\frac{\partial^1}{\partial q}+ 
\frac{\partial^1}{\partial q}\partial^0 = 0,\] and for any quandle 
we have:
\[\frac{\partial^0}{\partial q}\frac{\partial^1}{\partial q} + 
\frac{\partial^1}{\partial q}\frac{\partial^0}{\partial q} = 0.\]
\end{pf}

\begin{rem}
{\rm One could hope that
$\frac{\partial^i\partial^j }{\partial q\partial q'} + 
\frac{\partial^j\partial^i }{\partial q'\partial q} = 0$ for $i,j \in \{0,1\}$,
$q,q'\in Q$.  
This is not the case as the following example illustrates: \\
Let $w=(q_0,q_1,q_2,q_0,q_1,q_2)\in C^Q_6(R_3)$.
We have $\frac{\partial^1 w}{\partial q_1}=0$,
 while \[\frac{\partial^1 \partial^1 w}{\partial q_1 \partial q_2}= 
\frac{\partial^1 (q_1,q_0,q_2,q_1,q_0)}{\partial q_1} = 
-(q_0,q_2,q_1,q_0) \neq 0.\]}
\end{rem}

\subsection{Creating homology operations}

Let $w\in C_{\ell}^R(Q).$ We define a map 
$h_w\colon C_n^R(Q)\to C_{n+\ell}^R(Q)$ by:
\[h_w(x_1,\ldots,x_n)=(x_1,\ldots,x_n,w).\] 

Theorem \ref{main} below gives the necessary conditions on 
$w$, so that $(h_w)_*$ is a homology operation, in the language reminiscent of 
searching for extrema of a multi-variable function. Another criterion, 
applicable to non-connected quandles (e.g., $R_{2k}$), is given in 
Theorem \ref{Theorem 4.7}.

\begin{defn}
{\rm We say that a chain $w\in C_{\ell}^R(Q)$ is an {\it extreme chain} if $\partial^0(w)=0$, and all partial derivatives along $w$
are equal to zero ($\frac{\partial^1 w}{\partial q} = 0$, for all $q$).}
\end{defn}

\begin{thm}\label{main}
Assume that $w\in C_{\ell}^R(Q)$ is an extreme chain.  
Then,
$h_w$ is a chain map, and $(h_w)_*$ is a homological operation. 
Furthermore, because $h_w$ sends degenerate elements to degenerate elements,
$h_w$ is a chain map on quandle chains $C_n^Q(Q)$ ($\partial^0(w)$ and 
partial derivatives on quandle chains should be zero).
\end{thm}
\begin{pf} We have \[\partial h_w(u) = \partial (u,w) = 
(\partial u,w) +(-1)^{n+1}((u,\partial^0(w)) - 
(\sum_{q\in Q}(u*q,\frac{\partial^1 w}{\partial q})))=\]\[
h_w(\partial (u)) + (-1)^{n+1}((u,\partial^0(w))-
(\sum_{q\in Q}(u*q,\frac{\partial^1 w}{\partial q}))),\] and Theorem \ref{main}
follows. 

Notice that our proof can be interpreted as demonstrating 
that the map $h_{\partial^0 w}$ is chain homotopic to the map 
$u \to (\sum_{q\in Q}(u*q,\frac{\partial^1 w}{\partial q}))$.
\end{pf}

\begin{exmp}
{\rm Let $q_0,q_1,\ldots,q_{N-1}$ be a sequence of 
``Fibonacci elements" of $Q$, that is, $q_i*q_{i+1}=q_{i+2}$, with indices 
taken modulo $N$. Then:
\begin{enumerate}
\item[(i)]\cite{N-P-2} $s(q_0,q_1)=\sum_{i=0}^{N-1}(q_i,q_{i+1})$ satisfies 
the conditions of Theorem \ref{main}.
\item[(ii)] $w= \partial^0\frac{\partial^1 (q_0,q_1,...,q_{N-1},q_0)}
{\partial q_0}$ satisfies the conditions of Theorem \ref{main} for a quandle 
chain complex.
\end{enumerate}}
\end{exmp}

\begin{exmp}\ 
{\rm\begin{itemize}
\item[(i)]
The chain \[w=-(2,4,1)-(3,2,1)-(4,3,1)+(1,2,4)+(1,3,2)+(1,4,3),\]
where $1, 2, 3, 4$ denote the elements $0, 1, t$, and $1+t$ of the quandle $S_4=\Z_2[t]/(t^2+t+1)$, is an extreme chain, and for $g=(1,2)+(2,4)+(4,1)$, $(g,w)$ gives $\Z_2$ in $H_5^Q(S_4)$.
\item[(ii)]
For a dihedral quandle $R_3$, the chain \[w=-(1,0,1,2,0)-(1,2,0,2,0)-(2,0,2,1,0)-(2,1,0,1,0)\]\[+(0,1,0,1,2)+(0,1,2,0,2)+
(0,2,0,2,1)+(0,2,1,0,1)\] represents $\Z_3$ in $H_5^Q(R_3)$. It also satisfies the conditions of the Theorem \ref{main}, and the operation $h_w$ applied to low dimensional cycles gives torsion elements in homology groups of higher degrees.
\end{itemize}}
\end{exmp}

\begin{rem}\label{Remark 4.8}
{\rm In \cite{N-P-1} we introduced the concept of the $n$-th Burnside kei, that is a kei (involutive quandle) $Q$ in which
any pair of elements, $a_0,a_1$, satisfies the relation $a_n(a_0,a_1)=a_0$. Here we use the notation
\[a_n=a_n(a_0,a_1)=\underbrace{\ldots a_1*a_0*a_1*a_0*a_1}_{\mbox{total of n letters}}.\] 
Then the Fibonacci condition $a_{i+2}=a_i*a_{i+1}$ is satisfied, and we                              
can define the Fibonacci element $s(a_0,a_1)=\sum_{i=0}^{n-1}(a_i,a_{i+1})$ which is an extreme chain
in $C^R_2(Q)$, so it can be used to define a homology operation $h_{s(a_0,a_1)}$. In fact, 
the assumption that $Q$ is a kei is not needed here, so we define the $n$-th Burnside quandle to be the 
quandle for which the equation $a_n(a_0,a_1)=a_0$ holds for any pair of elements, $a_0$, $a_1$. 
In this setting, $s(a_0,a_1)$ is also an extreme chain and $h_{s(a_0,a_1)}$ is a homology operation; 
compare Example \ref{Example 5.6}.}
\end{rem}

\begin{thm}\label{Theorem 4.7}
 Let $Q$ be a quandle and $Q_1$ its invariant subquandle,
such that $*_q: Q \to Q$ does not depend on the choice of $q \in Q_1$.
Assume also that (any) $q\in Q_1$ is a $k$-element of $Q$ (i.e., $x*q^k=x,$ for any $x\in Q$). 
Then $h'_w = \sum_{i=0}^{k-1}h_w*q^i$ is a chain map for any 
$w \in {\Z}Q_1^{\ell}$.
\end{thm}
\begin{pf}
We have 
\[\partial h'_w(x_1,\ldots,x_n)=\partial \sum_{i=0}^{k-1}h_w*q^i = 
h'_w\partial(x_1,\ldots,x_n)\]\[ +(-1)^{n+1}\sum_{i=0}^{k-1}
((x_1,\ldots,x_n,\partial^0w)- ((x_1,\ldots,x_n)*q,\partial^1 w))*q^i.\]
Theorem \ref{Theorem 4.7} follows, as $Q_1$ is a trivial quandle with 
$\partial^0 w= \partial^1 w$.
\end{pf}

\begin{exmp}
{\rm The dihedral quandle $R_4$ is an example of a quandle which motivated 
Theorem \ref{Theorem 4.7} (here, $Q_1=\{0,2\}$). We generalize it as follows.
Consider two trivial quandles: $Q_0$ equal to $\Z_{k_0}$ as a set, and 
$Q_1$ equal to $\Z_{k_1}$ as a set. We define the quandle $Q_{k_0,k_1}$
on the set $Q_0\cup Q_1$ by $a*b= a+1$ for $a\in Q_0$ and $b \in Q_1$ or 
$a\in Q_1$ and $b \in Q_0$. Our quandle $Q$ (with $Q_1$ or $Q_0$ as invariant subquandle) satisfies 
the conditions of Theorem \ref{Theorem 4.7}. 
For example, we have three such quandles of size 6 ($Q_{5+1}$, $Q_{4+2}$, $Q_{3+3})$. Notice that $Q_{2,2}$ is the quandle $R_4$, and
the quandle $Q_{k_0,k_1}$ is an $lcm(k_0,k_1)$-quandle.}
\end{exmp}

\subsection{Naturality of homological operations}

The naturality we have in mind is defined as follows.
Suppose that we have a pre-homological operation for a given rack $X_1$, that is, 
a chain map $h\colon C^R_n(X_1) \to C^R_{n+k}(X_1)$, and there is a rack homomorphism $f\colon X_1 \to X_2$. 
Then, there is a uniquely defined pre-homological operation $T(h)\colon C^R_n(X_2) \to C^R_{n+k}(X_2)$, 
such that $f_{\#}h = T(h)f_{\#}$, where $f_{\#}\colon C^R_*(X_1) \to C^R_*(X_2)$ is the chain map 
induced by $f$ (we could also go straight to homology and write $ f_*h = (T(h))_*f_*$ as a 
condition for a homological operation).

Let $f\colon Q_1\to Q_2$ be a quandle homomorphism, and $f_*\colon H^W_n(Q_1)\to H^W_n(Q_2)$ be the induced homomorphism, where $W=R$, $Q$ or $D$. 
We define $f^{(n)}\colon Q_1^n\to Q_2^n$ by 
\[f^{(n)}(x_1,\ldots,x_n)=(f(x_1),\ldots,f(x_n)),\]
and extend this map linearly to the group homomorphism $f^{(n)}\colon {\Z}Q_1^n\to {\Z}Q_2^n$.\\
\begin{thm}\label{Theorem 4.10}\
\begin{enumerate}
\item[(i)] If $a$ is a $k$-element in $Q_1$, and $f(a)$ is a $k$-element in $Q_2$, then 
$fh_a'=h'_{f(a)}f,$ so $h_a'$ is natural.
\item[(ii)] If $w\in C_{\ell}^R(Q_1)$ is an extreme chain (as in Theorem \ref{main}), 
then \[(h_{f^{(\ell)}(w)})_*\colon H_n^R(Q_2)\to H_{n+\ell}^R(Q_2)\] 
is a homological operation.
\item[(iii)] The operation $h_w$ is natural.

\end{enumerate}
\end{thm}
\begin{pf}
(i) We have 
\[ f h'_a(x_1,\ldots,x_n)=
f(\sum_{i=0}^{k-1}(*_a)^i(x_1,\ldots,x_n,a))=\]\[\sum_{i=0}^{k-1}(*_{f(a)})^i(f(x_1),\ldots,f(x_k),f(a))=
h'_{f(a)}f,\] as needed.\\
(ii) 
Consider $\frac{\partial^1 f^{(\ell)}(w)}{\partial p}$. If $p$ is not in the image $f(Q_1),$ then $p$ is not an element occurring in $f^{(\ell)}(w)$, so 
$\frac{\partial^1 f^{(\ell)}(w)}{\partial p}=0$. Let $q_1,q_2,\ldots,q_k$ be the elements in the preimage $f^{-1}(p)$ that are in $w$.
We will show that \[ \frac{\partial^1 f^{(\ell)}(w)}{\partial p}=f^{(\ell-1)}\sum_{i=1}^{k} \frac{\partial^1 w}{\partial q_i}.\]
It is enough to show this formula in the case $w=(x_1,\ldots,x_{\ell}),$ where $x_i \in Q_1$ for $i=1,\ldots,\ell$:
\[f^{(\ell-1)}\sum_{i=1}^k \frac{\partial^1 (x_1,\ldots,x_\ell)}{\partial q_i}=\]
\[f^{(\ell-1)}\sum_{i=1}^k\sum_{j=1}^{\ell}(-1)^j
(x_1*q_i,\ldots,x_{j-1}*q_i,x_{j+1},\ldots,x_\ell)\delta_{x_j,q_i}=\]
\[\sum_{i=1}^k\sum_{j=1}^{\ell}(-1)^j f^{(\ell-1)}(x_1*q_i,\ldots,x_{j-1}*q_i,x_{j+1},\ldots,x_\ell)\delta_{x_j,q_i}=\]
\[\sum_{i=1}^k\sum_{j=1}^{\ell}(-1)^j(f(x_1)*f(q_i),\ldots,f(x_{j-1})*f(q_i),f(x_{j+1}),\ldots,f(x_\ell))\delta_{x_j,q_i}=\]
\[\sum_{j=1}^{\ell}(-1)^j(f(x_1)*p,\ldots,f(x_{j-1})*p,f(x_{j+1}),\ldots,f(x_\ell))\delta_{f(x_j),p}=\frac{\partial^1 f^{(\ell)}(w)}{\partial p}.\]
Similarly, if $\partial^0=0,$ then $\partial^0f^{(\ell)}(w).$\\
(iii) We have to check that $f^{(n+\ell)}h_w=h_{f^{(\ell)}(w)}f^{(n)}$. Indeed,
\[f^{(n+\ell)}h_w(x_1,\ldots,x_n)=f^{(n+\ell)}(x_1,\ldots,x_n,w)=(f(x_1),\ldots,f(x_n),f(w))=\]
\[h_{f^{(\ell)}(w)}(f(x_1),\ldots,f(x_n))=h_{f^{(\ell)}(w)}f^{(n)}(x_1,\ldots,x_n),\] as needed.
\end{pf}

\begin{rem}\label{Remark 4.12} 
{\rm Carter, Elhamdadi, and Saito \cite{CES} proposed a generalization of
 quandle homology to twisted quandle homology with $\partial^T = t\partial^0 -\partial^1$. 
Any extreme chain  $w\in C^R_{\ell}(Q)$ is also an extreme chain in the twisted homology (or any homology with the
differential of the form $a\partial^0 \pm b\partial^1$), and can be used to construct a natural
homology operation. This awaits detailed exploration.}
\end{rem}

\section{Further applications}\label{Section 5}
The techniques we have built so far can be used for some concrete calculations 
of quandle homology. We will illustrate it by several examples for both finite and infinite 
quandles. We start from a simple lemma crucial for our examples.
\begin{lem}\label{Lemma 5.1}
Consider a rack $X$ and its subrack $X_0$, with an embedding $i\colon X_0\to X$.
We say that a rack epimorphism $r\colon X\to X_0$ is a rack twist-retraction\footnote{Of course it is 
always possible, having a twist-retraction $r$, to consider a retraction $r'=(ri)^{-1}r$ 
(we have $r'i=(ri)^{-1}ri=Id_{X_0}$). However, it is sometimes convenient to work with twist-retractions of racks.} 
if $ri \in Aut(X_0)$, 
and a rack retraction if $ri = Id_{X_0}$. \\
Then $i_*\colon H^R_n(X_0) \to H^R_n(X)$ is a monomorphism, and $r_*\colon H^R_n(X) \to H^R_n(X_0)$ is 
an epimorphism.
\end{lem}

If $X$ is an Alexander quandle, then the retraction corresponds to an $\Z[t^{\pm 1}]$-module epimorphism 
which splits:
\begin{cor}\label{Corollary 5.2}
Consider a $\Z[t^{\pm 1}]$-module $A=A_1\oplus A_2$ with embeddings
$i_1\colon A_1\to A_1\oplus A_2$, $i_2\colon A_2\to A_1\oplus A_2$, and projections
$r_1\colon A_1\oplus A_2 \to A_1$, $r_2\colon A_1\oplus A_2 \to A_2$.
Naturally, $r_1i_1=Id=r_2i_2$, so $r_1$ and $r_2$ are retractions. We have:
\begin{enumerate}
\item[(i)] \[ (i_1)_*\colon H_n^W(A_1)\to H_n^W(A),\ (i_2)_*\colon H_n^W(A_2)\to H_n^W(A)\] are monomorphisms.
\item[(ii)] If $gcd(|A_1|,|A_2|)=1$, then
\[H_n^W(A_1)\oplus H_n^W(A_2)\] embeds in $H_n^W(A)$.
\end{enumerate}
\end{cor}
\begin{pf}
The first fact follows from Lemma \ref{Lemma 5.1}.
The second conclusion follows from the fact that $|A_i|^n$ annihilates $H_n^W(A_i)$, for $i=1, 2$.
\end{pf}

\begin{cor}\label{Corollary 5.3}\
\begin{enumerate}
\item[(i)] If $gcd(m,k)=1$, then $\Z_{mk}= \Z_{m}\oplus \Z_{k}$ and
\[ tor H_n^W(R_m)\oplus tor H_n^W(R_k)\subset H_n^W(R_{mk}).\]
\item[(ii)] For $k$ odd,
\[ tor H_n^W(R_k)\oplus tor H_n^W(R_k)\subset H_n^W(R_{2k}).\]
\item[(iii)] For $k$ odd, prime, such that $gcd(m,k)=1$,
$H^Q_3(R_{mk})$ and $H^Q_4(R_{mk})$, contain $\Z_k$ (compare Theorem \ref{Theorem 1.7} (ii)).
\end{enumerate}
\end{cor}

Inequality in Corollary \ref{Corollary 5.3}(ii) is seldom an equality, even for $k$ odd prime. For
example, $H^Q_4(R_6) = \Z_3^6= H^Q_4(R_3) \oplus H^Q_4(R_3)\oplus \Z_3^4$. The reason may be
that $R_6$, not being connected, allows for more homological operations than $R_3$ (e.g., two different $h'_a$
operations). We conjecture, however, the following:
\begin{conj}\label{Conjecture 5.4}\
\begin{enumerate}
\item[(i)]  \[tor (H^Q_{2n}(R_{4})) = (tor (H^Q_{2n-1}(R_{4})))^2 \oplus \Z_2^2 = \Z_2^{2(4^{n-1}-1)/3} \]
\[tor (H^Q_{2n+1}(R_{4})) = (tor (H^Q_{2n}(R_{4})))^2  = \Z_2^{4(4^{n-1}-1)/3}. \]
\item[(ii)] For $k$ odd:
\[tor (H^Q_{2n}(R_{2k})) = (tor (H^Q_{2n-1}(R_{2k})))^2 \oplus \Z_k^2 = \Z_k^{(5\cdot 4^{n-1}-2)/3},\ for\ n>1,\]
\[tor (H^Q_{2n+1}(R_{2k})) = (tor (H^Q_{2n}(R_{2k})))^2 =  \Z_k^{2 (5\cdot 4^{n-1}-2)/3}. \]
\item[(iii)]  \[2^{k-1}tor(H^Q_n(R_{2^k})) = \Z_2^{2f_n},\ \ \ for\ k>1.\]
Here, $\{f_n\}$ denotes the sequence of ``delayed" Fibonacci numbers, as in  Example 2.6.
\end{enumerate}
\end{conj}

\begin{exmp}\label{Example 5.5} 
{\rm Corollary \ref{Corollary 5.2} can be applied directly to 
an interesting family of Alexander quandles $A_{p,[k_1]_t[k_2]_t}= \Z_p[t^{\pm 1}]/([k_1]_t\cdot[k_2]_t)$,
with $gcd(k_1,k_2)=1$.  Recall the notation $[k]_t=1+t+\ldots+t^{k-1}$, $A_{p,[k]_t}=\Z_p[t^{\pm 1}]/([k]_t)$.\\
$A_{p,[k_1]_t[k_2]_t}$ splits to  $A_{p,[k_1]_t}\oplus  A_{p,[k_2]_t}$, so
\[tor H^W_n(A_{p,[k_i]_t})\subset torH_n^W(A_{p,[k_1]_t\cdot[k_2]_t}).\] For example, 
$A_{p,[2]_t}= R_p$, and $A_{2,[3]_t}=S_4$, so we know a lot about their quandle homology. 
We also computed, for instance, that $H^Q_2(A_{2,[5]_t})=\Z_2^2 $, 
and $H^Q_2(A_{2,[4]_t})= H^Q_2(A_{2,[6]_t}) = \Z^2 \oplus \Z_2^4$.}
\end{exmp}

\begin{figure}
\begin{center}
\includegraphics[height=3.3 cm]{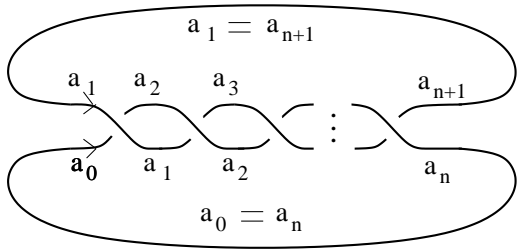}
\caption{The torus link $T_{2,n}$ with Fibonacci quandle labeling}\label{torus link}
\end{center}
\end{figure}

\begin{exmp}\label{Example 5.6}
{\rm Consider the Fibonacci quandle $F_n=\{a_0,a_1,\ldots,a_{n-1}\ | \ a_2=a_0*a_1, a_3=a_2*a_1,\ldots,a_{n-1}=a_{n-2}*a_{n-3},
a_0=a_{n-2}*a_{n-1},a_1=a_{n-1}*a_0\}$. Then we have a Fibonacci element $s(a_0,a_1)=\sum_{i=0}^{n-1} (a_i,a_{i+1})$
which can be used for a homological operation $h_{s(a_0,a_1)}$. We also have an epimorphism
$r\colon F_n \to R_n$ given by $r(a_j)=j$. This can be used to produce nontrivial elements in
$H^W_k(F_n)$. Notice that $F_n$ is the fundamental quandle of the torus link of type $(2,n)$, $T_{2,n}$;
compare Figure \ref{torus link}.
\begin{enumerate}
\item[(i)] $H^Q_n(F_4)$ contains the free part for $2\leq n \leq 7$, because we found that
in $H^Q_n(R_4)$ the elements:\\
$s(0,1),\ h_{s(0,1)}((0)),\ h_{s(0,1)}(s(0,1)),\ h_{s(0,1)}h_{s(0,1)}((0))$,\\
$h_{s(0,1)}h_{s(0,1)}(s(0,1)),$ and $h_{s(0,1)}h_{s(0,1)}h_{s(0,1)}((0))$ represent $\Z$. Compare \cite{Eis}, 
where it was shown that $H^Q_2(F_4)=\Z^2$.
\item[(ii)] $H_3^Q(F_p)$ is nontrivial for $p$ odd prime. It follows from the fact that\\
$r_*((a_0,s(a_0,a_1)))$ is a non-zero element in $H^Q_3(R_p)=\Z_p$. In fact, in all cases we checked
$h_{s(0,1)}\colon C^Q_n(R_k) \to C^Q_{n+2}(R_k)$ induces a monomorphism on homology \cite{N-P-2}.
From this follows that $H^Q_{2n+1}(F_3)$ is nontrivial for $n\leq 5$. One can conjecture that
$h_{s(a_0,a_1)}\colon C^Q_n(F_k) \to C^Q_{n+2}(F_k)$ induces a monomorphism on homology. Our methods
are not yet sufficient to decide whether $H^Q_n(F_k)$ can have torsion elements. To answer this question, we plan to
use the detailed description of $F_3$ in \cite{N-P-3} (the fundamental quandle of
the trefoil knot was interpreted there as a symplectic quandle).
\end{enumerate}}
\end{exmp}

We end with a few rather general conjectures (\ref{Conjecture 5.7}, \ref{Conjecture 5.9}, \ref{Conjecture 5.10}), 
and two new examples illustrating them:
\begin{conj}\label{Conjecture 5.7}
Let $r\colon Q_1 \to Q_2$ be an epimorphism of finite connected quandles, such that
$|Q_1|= k|Q_2|$, and $gcd(k,|Q_2|)=1$. Then,
$r_*\colon H^W_n(Q_1) \to H^W_n(Q_2)$ is an epimorphism, and if restricted to
the part of $H^R_n(Q_1)$ not annihilated by $k^j$ for some $j$, it is an isomorphism.
\end{conj}

\begin{exmp}\label{Example 5.8} 
{\rm Consider the 6-element connected quandles $Q_{1}$ and $Q_2$ which are the last two quandles 
from the page 195 of \cite{CKS}. We have an epimorphism $r_1\colon Q_1 \to R_3$, and $r_2\colon Q_2 \to R_3$.
As far as computation of homology was performed, Conjecture \ref{Conjecture 5.7} holds. In particular,
non-2 torsion part of homology of $Q_i$, $i=1,2$, is the same as that for $R_3$.\\
For example, $Q_2$, the last quandle in \cite{CKS}, is a 4-quandle, given by the following $*$-multiplication 
table:
\begin{center}
\begin{tabular}{c| c ccccc} 
$*$&1 & 2 & 3 & 4 & 5 & 6  \\
\hline 
1 & 1 & 1 & 6 & 3 & 4 & 5 \\
2 & 2 & 2 & 4 & 5 & 6 & 3 \\
3 & 4 & 6 & 3 & 3 & 3 & 1 \\
4 & 5 & 3 & 1 & 4 & 2 & 4 \\
5 & 6 & 4 & 5 & 1 & 5 & 2 \\
6 & 3 & 5 & 2 & 6 & 1 & 6 \\
\end{tabular}
\end{center} 
\noindent and $r_2\colon Q_2 \to R_3$
is given by $r_2(1)=r_2(2)=0$, $r_2(3)=r_2(5)=1$, and $r_2(4)=r_2(6)=2$.
$s(1,3)= (1,3)+(3,6)+(6,1)$ is a Fibonacci chain in $C_2^Q(Q_2)$, with 
$r_2(s(1,3))= s(0,1)\in C_2^Q(R_3)$.}
\end{exmp}
If we combine Example \ref{Example 5.8} with the prediction that $3\cdot 2^i$ annihilates torsion in homology 
of $Q_1$ and $Q_2$, then we can suggest:

\begin{conj}\label{Conjecture 5.9}
Let oddtor denote the odd part of the torsion of a group. Then
\[oddtor(H^W_n(Q_1)) = oddtor(H^W_n(Q_2)) = oddtor(H^W_n(R_3)) = \Z_3^{f_n},\]
where $f_n$ is as defined before.
\end{conj}

The last, rather general conjecture should be the first step to understand torsion 
in homology of quandles, but, until now, the only nontrivial quandle for which it is solved 
is the dihedral quandle $R_3$ \cite{N-P-2}.
\begin{conj}\label{Conjecture 5.10}
If $Q$ is a finite quasigroup quandle, then $|Q|$ annihilates torsion of its homology.
\end{conj}
Notice that {\it quasigroup} cannot be replaced by {\it connected} as 6 does not annihilate the torsion of homology of $Q_2$. In particular, $H^Q_3(Q_2)=\Z_8\oplus\Z_3$ (see \cite{CKS-1}).

\section{Acknowledgements}
M.~Niebrzydowski was partially supported by the Louisiana Board of Regents grant (\# LEQSF(2008-11)-RD-A-30), 
and by the research award from the University of Louisiana at Lafayette.

J.~H.~Przytycki was partially supported by the NSA grant (\# H98230-08-1-0033), 
by the NSF grant (\# 0745204), and by the CCAS/UFF award.

We would like to thank the referee  and T.~Nosaka for many suggestions that improved our paper.


\begin{thebibliography}{99}

\bibitem{A-G} N.~ Andruskiewitsch, M.~Gra\~{n}a, 
From racks to pointed Hopf algebras, Adv. Math. 178 (2003) 177-243.

\bibitem{CES}
S.~Carter, M.~Elhamdadi, M.~Saito, Twisted quandle homology theory and cocycle knot invariants, 
Algebraic and Geometric Topology 2 (2002) 95-135.

\bibitem{CJKLS}
S.~Carter, D.~Jelsovsky, S.~Kamada, L.~Langford, M.~Saito,
State-sum invariants of knotted curves and surfaces from quandle cohomology,
Electron. Res. Announc. Amer. Math. Soc. 5 (1999) 146-156.

\bibitem{CJKS-1}
S.~Carter, D.~Jelsovsky, S.~Kamada, M.~Saito,
Quandle homology groups, their Betti numbers, and virtual knots,
J. Pure Appl. Algebra 157 (2001) 135-155.


\bibitem{CJKS-2}
S.~Carter, D.~Jelsovsky, S.~Kamada, M.~Saito,
Shifting homomorphisms in quandle cohomology and skeins of cocycle knot invariants,
Journal of Knot Theory and its Ramifications 10 (2001) 579-596.

\bibitem{CKS}
S.~Carter, S.~Kamada, M.~Saito, Surfaces in 4-space, Encyclopaedia
of Mathematical Sciences, Low-Dimensional Topology III, R.V.Gamkrelidze,
V.A.Vassiliev, Eds., Springer-Verlag Berlin-Heidelberg, 2006, 213pp.

\bibitem{CKS-1}
S.~Carter, S.~Kamada, M.~Saito, Geometric interpretations of quandle homology, 
Journal of Knot Theory and its Ramifications 10 (2001) 345-386.

\bibitem{E-G}
P.~Etingof, M.~Gra\~{n}a, On rack cohomology,
J. Pure Appl. Algebra 177 (2003) 49-59.


\bibitem{Eis}
M.~Eisermann, Homological characterization of the unknot, 
J. Pure Appl. Algebra 177 (2003) 131-157.

\bibitem{F-R}
R.~Fenn, C.~Rourke,
Racks and links in codimension two,
Journal of Knot Theory and its Ramifications 1 (1992) 343-406.

\bibitem{FRS}
R.~Fenn, C.~Rourke and B.~Sanderson, James bundles and applications,
preprint.\\
e-print:\ {\tt http://www.maths.warwick.ac.uk/$\sim$cpr/ftp/james.ps}


\bibitem{GAP}
The GAP Group, GAP -- Groups, Algorithms, and Programming,\\
http://www.gap-system.org.

\bibitem{Gr} M. Greene, Some results in geometric topology and geometry, 
Ph.D. thesis, University of Warwick, advisor: Brian Sanderson, 1997.

\bibitem{JoyTh} D.~Joyce, An algebraic approach to symmetry with applications to knot theory, Ph.D. thesis, University of Pennsylvania, advisor: Peter Freyd, 1979.

\bibitem{Joy} D.~Joyce, A classifying invariant of knots: the knot
quandle, J. Pure Appl. Algebra 23 (1982) 37-65.

\bibitem{Kam} S.~Kamada,
Quandles with good involutions, their homologies and knot invariants, in: Intelligence of Low Dimensional Topology 2006, Eds. J. S. Carter et. al., World Scientific Publishing Co., 2007, pp. 101-108.

\bibitem{L-N}
R.~A.~Litherland, S.~Nelson, The Betti numbers of some finite racks,
J. Pure Appl. Algebra 178 (2003) 187-202.

\bibitem{Mat}
S.~Matveev, Distributive groupoids in knot theory, 
Math. Sbornik 119 (1982) 78-88. Eng. transl.:  
Math. USSR Sbornik 47 (1984) 73-83.

\bibitem{Moc}
T.~Mochizuki, 
Some calculations of cohomology groups of finite Alexander quandles,
J. Pure Appl. Algebra 179 (2003) 287-330.


\bibitem {N-P-1}
M.~Niebrzydowski, J.~H.~Przytycki, Burnside Kei,
Fundamenta Mathematicae 190 (2006) 211-229.



\bibitem {N-P-2}
M.~Niebrzydowski, J.~H.~Przytycki, 
Homology of dihedral quandles,
J. Pure Appl. Algebra 213 (2009) 742-755.


\bibitem {N-P-3}
M.~Niebrzydowski, J.~H.~Przytycki, 
The quandle of the trefoil as the Dehn quandle of the torus,
Osaka Journal of Mathematics 46 (2009) 645-659.


\bibitem{Oht}
T.~Ohtsuki, Quandles, in: Problems on invariants of knots and 3-manifolds,
Geometry and Topology Monographs 4 (2003) 455-465.


\bibitem{Tak}
M.~Takasaki, Abstraction of symmetric transformation, (in Japanese)
Tohoku Math. J. 49 (1942/3) 145-207.

\end{thebibliography}
\end{document}